\documentclass[12pt]{article}
\usepackage{amssymb, amsmath, latexsym, a4} 
\usepackage[latin1]{inputenc}
\usepackage[all]{xy}

\begin{document}

\newcommand{\qed}{\hfill $\Box $}
\def \G{\mathop{\Gamma }\nolimits}
\def \Y{\mathop{\cal Y }\nolimits}
\def \n{\mathop{\underline n }\nolimits}
\def \m{\mathop{\underline m }\nolimits}
\newcommand{\Ga}{\raisebox{0.1mm}{$\Gamma$}}
\def\S{\smallskip \par}
\def\M{\medskip \par}
\def\ll{ \lambda}
\def\B{\bigskip \par}
\def\BB{\bigskip \bigskip \par}
\newcommand{\Q}{\mathbb{Q}}
\newcommand{\F}{\mathbb{F}}

\newcommand{\Z}{\mathbb{Z}}
\def\t{\otimes }
\def\ee{\epsilon }
\newcommand{\pn}{\par \noindent}
\def \H{\mathop{\sf H}\nolimits}
\def \Im{\mathop{\sf Im}\nolimits}
\def \Ker{\mathop{\sf Ker}\nolimits}
\def \Coker{\mathop{\sf Coker}\nolimits}
\def \Hom{\mathop{\sf Hom}\nolimits}
\def \Ext{\mathop{\sf Ext}\nolimits}
\def \Tor{\mathop{\sf Tor}\nolimits}

\newtheorem{De}{Definition}[section]
\newtheorem{Th}[De]{Theorem}
\newtheorem{Pro}[De]{Proposition}
\newtheorem{Le}[De]{Lemma}
\newtheorem{Co}[De]{Corollary}
\newtheorem{Rem}[De]{Remark}
\newtheorem{Ex}[De]{Example}
\newcommand{\Def}[1]{\begin{De}#1\end{De}}
\newcommand{\Thm}[1]{\begin{Th}#1\end{Th}}
\newcommand{\Prop}[1]{\begin{Pro}#1\end{Pro}}

\newcommand{\Lem}[1]{\begin{Le}#1\end{Le}}
\newcommand{\Cor}[1]{\begin{Co}#1\end{Co}}
\newcommand{\Rek}[1]{\begin{Rem}#1\end{Rem}}

\BB
\centerline{\large{\bf Andr\'e-Quillen  homology via functor homology}}
\bigskip

\centerline{\bf Teimuraz Pirashvili\footnote{The author was partially supported by the grant
INTAS-99-00817 and by the \pn TMR network K-theory and
algebraic groups , ERB FMRX CT-97-0107.} }

\BB
{\footnotesize We obtain  Andr\'e-Quillen homology for commutative
algebras using relative homological algebra in the category of functors
on finite pointed sets.}

\bigskip

\pn {\bf Keywords:}
Andr\'e-Quillen  homology, relative homological algebra, functor homology, 
Gamma spaces \pn

\pn {\bf AMS classification:} 13D03, 18G25.

\BB

 \section { Introduction.} 
\pn Let $\Ga$ be the
small category of
finite pointed sets. For any $n\geq 0$, let $[n]$ be the set
$\lbrace 0,1,...,n\rbrace$ with basepoint $0$. We  assume that  the 
objects of $\Ga$ are the sets $[n]$. A left $\Ga$-module is  a covariant 
functor $\Ga \to \sf Vect$ to the category of vector spaces over a field
$K$.  For a left
$\Ga$-module $F$ we put 
$$\pi_0(F):=\Coker(d_0-d_1+d_2:F([2])\to F([1])),$$
where $d_1$ is induced by the folding map $[2]\to [1]$, $1,2\mapsto 1$
while $d_0$ and $d_2$ are induced by the projection maps 
$[2]\to [1]$ given respectively
by $1\mapsto 1, 2\mapsto 0$ and $1\mapsto 0, 2\mapsto 1$.
The category $\Ga$-$mod$ of left $\Ga$-modules is an abelian
category with enough projective and injective objects. Therefore one can 
form the left derived functors of the functor $\pi_0:\Ga$-$mod\to \sf Vect$,
which we will denote by $\pi_*$. Thanks to \cite{hodge} and \cite{dold}
we know that $\pi_*F$ is isomorphic to the homotopy  of the spectrum
corresponding to the $\Ga$-space $F$ according to Segal (see \cite{segal}
and \cite{BF}). 

Let $A$ be a commutative algebra
over a ground field $K$ and let $M$ be an $A$-module. There exists a 
functor ${\cal L}(A,M):\Ga\to {\sf Vect}$,
which assigns $M\t A^{\t n}$ to $[n]$ (see  \cite{JLL} or section 3). 
Here  all tensor products are  
taken over $K$. 
It was proved in \cite{piri} that $\pi_*({\cal L}(A,M))$ is 
isomorphic to a brave new algebra version
of Andr\'e-Quillen homology ${\sf H}_*^{\G}(A,M)$ constructed by 
Alan Robinson and
Sarah Whitehouse \cite{Sarah}. The main result of this paper shows 
that a similar
isomorphism exists also for 
Andr\'e-Quillen homology if one takes an appropriate
relative derived functors of the same functor $\pi_0:\Ga$-$mod\to \sf Vect$.

\eject
\BB
\section{ A class of proper exact sequences}
\pn Thanks to
the Yoneda lemma  $\Ga^n$, $n\geq 0$, 
are projective generators of
the category $\Ga$-$mod$. Here 
$$\Ga ^n:\ =K[\Hom_{\G}([n],-)] .$$
and $K[S]$ denotes the free vector space  generated by a set $S$.
For left 
$\Ga$-modules $F$ and $T$ one defines the pointwise tensor product
 $F\t T$  to be the left $\Ga$-module given 
by $(F\t T)([n])= F([n])\t T([n])$. Since
$\Ga ^n\t \Ga^m\cong \Ga ^{n+m}$ one sees that the 
tensor product of two projective
left $\Ga$-modules is still projective. We also have
$\Ga^n\cong (\Ga^1)^{\t n}.$

\M
{\it A partition} $\ll =
(\ll_1,\cdots, \ll_k)$
is a sequence
of natural numbers $\ll_1\geq \cdots  \geq \ll_k\geq 1$. 
The sum of partition is given by
$s(\ll):=\ll_1+\cdots +\ll_k$, while the group $\Sigma (\ll)$ 
is a product of the corresponding symmetric groups 
$$\Sigma ({\ll}):=\Sigma_{\ll_1}\times \cdots \times \Sigma_{\ll_k}.$$
which is identified with the Young subgroup of $\Sigma _{s(\ll)}$. 
Let us observe that $\Sigma _n=Aut_{\G}([n])$ and therefore $\Sigma _n$
acts on $\Ga^n\cong (\Ga^1)^{\t n}$. For a 
partition $\ll$ with $s(\ll)=n$ we let $\Ga ({\ll})$ be the 
coinvariants of $\Ga^n$ under the action of 
$\Sigma(\ll)\subset \Sigma _n$. 

For a vector space $V$
we let $S^*(V)$, $\Lambda ^*(V)$ and $D^*(V)$ be respectively 
the symmetric, exterior and 
divided power algebra generated by $V$. Let us recall that 
$S^n(V)=(V^{\t n})/\Sigma _n$ is the space of coinvariants of $V^{\t n}$
under the
action of the symmetric group $\Sigma _n$, while
 $D^n(V)=(V^{\t n})^{\Sigma _n}$ is the space of invariants. Moreover for a 
particiaon $\ll=(\ll_1,\cdots , \ll_k)$ we put 
$$S^{\ll}:=S^{\ll_1}\t \cdots \t S^{\ll_k}.$$ A similar meaning has also
$\Lambda ^{\ll}$ and $D^{\ll}$. It follows from the definition that
$$\Ga(\ll)\cong S^{\ll}\circ \Ga^1.$$
In particular $\Ga(1,\cdots, 1)\cong \Ga^n$ and $\Ga(n)\cong S^n\circ \Ga^1.$

\M
Let $$0\to T_1\to T\to T_2\to 0 $$
be an exact sequence of left $\Ga$-modules. It is called a {\it $\Y$-exact 
sequence} if for any partition $\ll$ with $s(\ll)=n$ the 
induced map 
$$T([n])^{\Sigma (\ll)}\to T_2([n])^{\Sigma (\ll)}$$ 
is surjective. Here and elsewhere, $M^G$ denotes 
the subspace of $G$-fixed elements of a $G$-module $M$.
For a $\Y$-exact sequence $0\to T_1\to T\to T_2\to 0$  the sequence
$$0\to T_1([n])^{\Sigma (\ll)}\to T([n])^{\Sigma (\ll)}\to T_2([n])
^{\Sigma (\ll)}\to 0 $$
is also exact. Following to Section XII.4 of 
\cite{homology} we introduce the related
notions. 
An epimorphism $f:F\to T$ is called {\it $\Y$-epimorphism}
if $$0\to \Ker (f)\to F\to T\to 0$$ is a $\Y$-exact sequence. Similarly, 
a monomorphism $f:F\to T$
 is called {\it $\Y$-monomorphism}
if $$0\to  F\to T\to \Coker (f)\to 0$$ is a $\Y$-exact sequence. 
A morphism $f:F\to T$ is called  {\it $\Y$-morphism} if $F\to \Im(f)$ is
a $\Y$-epimorphism and $\Im (f)\to T$ is a $\Y$-monomorphism.
A left $\Ga$-module $Z$ is called {\it $\Y$-projective} if for any 
$\Y$-epimorphism $f:F\to T$ and any morphism $g:Z\to T$ there exist
a morphism $h:Z\to F$ such that $g=fh$.

\begin{Le} {\rm i)} If a short exact sequence is isomorphic to
a $\Y$-exact sequence, then it is also a $\Y$-exact sequence.

\S
{\rm ii)} A split short exact sequence is  $\Y$-exact.

\S
{\rm iii)} A composition of two $\Y$-epimorphisms is still a $\Y$-epimorphism.

\S
{\rm iv)} If $f$ and $g$ are two composable epimorphism and $fg$ is a 
$\Y$-epimorphism, then $f$ is also a $\Y$-epimorphism.

\S
{\rm v)} A composition of two $\Y$-monomorphisms is still a $\Y$-monomorphism.

\S
{\rm vi)} If  $f$ and $g$ are two composable monomorphism and $fg$ is a 
$\Y$-monomorphism, then $g$ is also a $\Y$-monomorphism.

\end{Le}

\S
\pn {\it Proof}. The properties i)- iv) are clear.
Let $f:B\to C$ and $g:A\to B$ be monomorphisms. One can form the 
following diagram

$$  \xymatrix{& & 0 \ar[d] & 0\ar[d] \\ 
0 \ar[r]& A \ar[r]^{g}\ar[d]^{1_A}& B\ar[r]\ar[d]^{f} &X\ar[r]\ar[d]& 0\\
0 \ar[r]& A \ar[r]^{fg}& C\ar[r]\ar[d] &Z\ar[r]\ar[d]& 0\\
&& Y\ar[r]^{1_Y}\ar[d]& Y\ar[d]\\
&& 0 &0\\
}$$
Assume $f$ and $g$ are are $\Y$-monomorphisms, then for any partition
$\ll$ with $s(\ll)=n$ one has a commutative diagram
$$  \xymatrix{& & 0 \ar[d] & 0\ar[d] \\ 
0 \ar[r]& A([n])^{\Sigma (\ll)} \ar[r]\ar[d]^{1_A}& B([n])^{\Sigma (\ll)}\ar[r]\ar[d] &X([n])^{\Sigma (\ll)}\ar[r]\ar[d]& 0\\
0 \ar[r]& A ([n])^{\Sigma (\ll)}\ar[r]& C([n])^{\Sigma (\ll)}\ar[r]^{h}\ar[d] &Z([n])^{\Sigma (\ll)}\ar[d]\\
&& Y([n])^{\Sigma (\ll)}\ar[r]^{1_Y}\ar[d]& Y([n])^{\Sigma (\ll)}\\
&& 0 \\
}$$
The diagram chasing shows that $h$ is an epimorphism and therefore $fg$ is a $\Y$-monomorphism and v) is proved. Assume now $fg$ is a $\Y$-monomorphism. Them we have the 
following commutative diagram
$$  \xymatrix{& & 0 \ar[d] & 0\ar[d] \\ 
0 \ar[r]& A([n])^{\Sigma (\ll)} \ar[r]\ar[d]^{1_A}& B([n])^{\Sigma (\ll)}
\ar[r]^{l}\ar[d] &X([n])^{\Sigma (\ll)}\ar[d]\\
0 \ar[r]& A ([n])^{\Sigma (\ll)}\ar[r]& C([n])^{\Sigma (\ll)}\ar[r]\ar[d] &Z([n])^{\Sigma (\ll)}\ar[d]\ar[r]&0\\
&& Y([n])^{\Sigma (\ll)}\ar[r]^{1_Y}& Y([n])^{\Sigma (\ll)}\\
 }$$
The diagram chasing shows that $l$ is an epimorphism and therefore $f$ is a $\Y$-monomorphism and therefore we get vi). \qed

\S
As an immediate corollary we 
obtain that the class of all ${\Y}$-exact sequences
is proper in the sense of Mac Lane \cite{homology}. We now show that
there are enough $\Y$-projective objects.

\S

\begin{Le}\label{proj} {\rm i)} For any partition 
$\ll$ the left $\Ga$-module $\Ga({\ll})$
is a $\Y$-projective object.

\S
{\rm ii)} A morphism $f:F\to T$ of left $\Ga$-modules is $\Y$-epimorphism
iff for any partition $\ll$ the induced morphism
$$\Hom_{\G-mod}(\Ga({\ll}),F)\to \Hom_{\G-mod}(\Ga({\ll}),T)$$
is an epimorphism.

\S
{\rm iii)} For any left $\Ga$-module $F$ there is a $\Y$-projective
object $Z$ and $\Y$-epimorphism $f:Z\to F$.

\S
{\rm iv)} Any projective $\Y$-module is a direct summand of the sum of
objects of the form  $\Ga({\ll})$.

\S
{\rm v)} The tensor product of two $\Y$-projective 
left $\Ga$-modules is still $\Y$-projective.

\end{Le}

\M
\pn {\it Proof}. Let $\ll$ be a partition with $s(\ll)=n$. By definition
$\Ga({\ll})=H_0(\Sigma ({\ll}),\Ga ^n)$. Hence for any left $\Ga$-module $F$
one has
$$\Hom_{\G-mod}(\Ga({\ll}),F)\cong H^0(\Sigma({\ll}), \Hom_{\G-mod}(\Ga^n,F)
\cong F(n)^{\Sigma(\ll)}.$$ 
The assertions i) and ii) are immediate consequence of this isomorphisms. 
To proof iii) we set
$$X(\ll):=\Hom_{\G-mod}(\Ga({\ll}), F).$$
Moreover, for each $x\in X(\ll)$ we let $f_x:\Ga({\ll})\to F$ be the
corresponding morphism. Take 
$Z=\bigoplus_{\ll}\bigoplus_{x\in X(\ll)}\Ga ({\ll})$. 
Then the collection $f_x$,
$x\in X(\ll)$ 
yields the morphism $f:Z\to F$. We have to
show that it is  a $\Y$-epimorphism. Let $g:\Ga({\ll})\to F$
be a morphism of left $\Ga$-modules. By ii) we need to lift $g$ to $Z$.
By our construction $g\in X(\ll)$ and therefore the inclusion
 $\Ga({\ll})\to Z$
corresponding to the summand $g\in X(\ll)$ 
is an expected lifting and iii) is proved. 
The proof of iii) shows that one can assume
$P$ to be a sum of $\Ga^{\ll}$ and iv) follows. To proof the last 
statement one
observes that, for any partitions  $\ll$ and $\mu$ 
one has
$$\Ga ({\ll})\t \Ga ({\mu})\cong (\Ga^{s(\ll)}\t 
\Ga^{s(\mu)})^{\Sigma (\ll)\times
\Sigma (\mu)}=(\Ga^{s(\ll)+s(\mu)})^{\Sigma (\ll)\times
\Sigma (\mu)}$$
and therefore $\Ga ({\ll})\t \Ga ({\mu})$ is $\Y$-projective. \qed

 


\section{ Definition of Andr\'e-Quillen homology and the functor $\cal L$}

\pn The definition  Andr\'e-Quillen
homology is based on the framework of homotopical algebra \cite{Q}
and it is given as follows. We let $C_*(V_*)$ be the
chain complex associated to a simplicial vector space $V_*$. 
Let $A$ be a commutative algebra
over a ground field $K$ and let $M$ be an $A$-module. {\it 
A simplicial resolution of $A$} is an augmented  simplicial object
$P_*\to A$ in the category of 
commutative algebras, which is a weak equivalence (in other words
$C_*(P_*)\to A$ is a weak equivalence). A simplicial resolution
is called {\it free} if   
$P_n$ is a polynomial algebra over $K$
for all $n\geq 0$. Any commutative algebra posses a free simplicial
resolution which is unique up to homotopy. 
Then the Andr\'e-Quillen
homology is defined by
$${\sf D}_*(A,M):=H_*(C_*(\Omega^1_{P_*}\t_{P_*}M)),$$
where $\Omega^1$ is the K\"ahler 1-differentials and $P_*\to A$ is
a free simplicial resolution. In the dimension
$0$ we have ${\sf D}_0(A,M)\cong \Omega^1_{A}\t_{A}M.$

\M
As we mentioned the functor  ${\cal L}(A,M):\Ga\to \sf Vect$ 
is given on objects by $[n]\mapsto M\t A^{\t n}.$
For a pointed map $f:[n]\to [m]$, the action of 
$f$ on ${\cal L}(A,M)$ is given by
$$f_*(a_0\t \cdots \t a_n):\ = b_0\t \cdots \t b_m,$$
where
$b_j=\prod _{f(i)=j}a_i, \ j=0,\cdots ,n. $

\begin{Ex}\label{patara} 
{\rm Let $M=A=K[t]$. In this case one has an isomorphism
$${\cal L}(K[t],K[t])\cong S^*\circ \Ga^1.$$
To see this isomorphism, one observes that $\Ga^1$ assigns the free vector 
space on a set $[n]$ to $[n]$ and therefore both functors in the question
assign the ring $K[t_0,\cdots , t_n]$ to $[n]$. An important consequence
of this isomorphism is the fact that the functor 
${\cal L}(K[t],K[t])$ is $\Y$-projective.} 
\end{Ex}

\begin{Le}\label{nuli}
 For any commutative algebra $A$ and any $A$-module $M$, 
one has a natural isomorphism
$\pi_0({\cal L}(A,M))\cong \Omega^1_{A}\t_{A}M.$
\end{Le}
\S
\pn {\it Proof}. This is a consequence of Proposition 1.15 and Proposition 2.2 
of \cite{hodge}. \qed

\begin{Le}\label{sizuste}
{\rm i)} Let $A$ be a commutative algebras and let
$$0\to M_1\to M\to M_2\to 0$$
be a short exact sequence of $A$-modules.
Then
$$0\to {\cal L}(A,M_1)\to {\cal L}(A,M)\to {\cal L}(A,M_2)\to 0$$
is a $\Y$-exact sequence.
\S
{\rm ii)} Let $f:B\to A$ be a surjective homomorphism of commutative algebras, then for any
$A$-module $M$ the induced morphism of left $\Ga$-modules 
$${\cal L}(B,M)\to 
{\cal L}(A,M)$$ is a $\Y$-epimorphism.
\end{Le}
\S
\pn {\it Proof}.  One observes that for any partition $\ll$ with $s(\ll)=n$
one has
$$({\cal L}(A,M)([n]))^{\Sigma (\ll)}= (M\t A^{\t n})^{\Sigma (\ll)}
\cong M\t D^{\ll}(A).$$
Since we are over field the tensor product is exact and we obtain i). 
By the same reason $f$  has a linear section, which yields
also a linear section of  $D^{\ll}(B)\to D^{\ll}(A)$, because $D^{\ll}$ is a functor
defined on the category of vector spaces .\qed
 
\BB
\section { Relative derived functors}
\pn By Lemma \ref{proj} the
class of $\Y$-exact sequences has enough projective objects. Thanks to 
\cite{homology} this allows us to construct the relative derived functors. 
Let us recall that an augmented chain complex $X_*\to F$ is called 
a {\it $\Y$-resolution}
of $F$ if it is exact (that is $H_i(X_*)=0$, for $i>0$
and $H_0(X_*)\cong F$) and all boundary maps $X_{n+1}\to X_n$ are
$\Y$-morphisms, $n\geq 0$. It follows from  Lemma \ref{proj} that $X_*\to F$
is a $\Y$-resolution iff for any partition $\ll$ the
augmented complex 
$$\Hom_{\G-mod}(\G({\ll}),X_*)\to \Hom_{\G-mod}(\G({\ll}),F)$$
is exact. A $\Y$-resolution $Z_*\to F$ is called {\it $\Y$-projective 
resolution} if for all $n\geq 0$ the left $\Ga$-nodule
 $Z_n$ is a $\Y$-projective object. We define 
$\pi^{\Y}_*(F)$ using relative derived functors of the functor $\pi_0:\Ga$-$mod\to \sf Vect$. In other words we 
put
$$\pi^{\Y}_n(F):=H_n(\pi_0(Z_*)), \ n\geq 0,$$
where $Z_*\to F$ is a $\Y$-projective
resolution. By \cite{homology} this gives the well-defined functors
$\pi_n^{\Y}:\Ga$-$mod\to \sf Vect$, \ $n\geq 0$.

\Lem{If $K$ is a field of chracteristic zero, then 
$\pi_*(F)\cong \pi_*^{\Y}(F)$.}

\S
\pn{\it Proof}. In this case all exact sequences are $\Y$-exact, because
for any finite group $G$, the functor $M\mapsto M^G$ is exact.\qed



\begin{Le}\label{namravli}For left $\Ga$-modules $F,T$ one has  
an isomorphism
$$\pi^{\Y} _*(F\t T)\cong \pi^{\Y}_*(F)\t T([0]) 
\oplus F([0])\t \pi^{Y} _*(T).$$  
\end{Le}

\S
\noindent {\it Proof}. The result in the dimension $0$ is known (see Lemma 4.2
of \cite{hodge}). Let $Z_*\to F$
and $R_*\to T$ be $\Y$-projective resolutions. By Lemma \ref{proj}
 $Z_* \t R_*\to F\t T$
is aslo $\Y$-projective resolution. Thus 
$$\pi^{\Y} _*(F\t T)=H_*(\pi_0(Z_*\t R_*))\cong$$
$$H_*( 
\pi^{\Y}_0(Z_*)\t R_*([0]) \oplus Z_*([0])\t \pi^{Y} _0(R_*))\cong$$
$$\pi^{\Y}_*(F)\t T([0]) \oplus F([0])\t \pi^{Y} _*(T),$$
where the last isomorphism follows from the
Eilenberg-Zilber theorem and K\"unneth theorem. \qed

\begin{Le}\label{rezolventa}
Let $\ee:X_*\to A$ be a simplicial resolution in the category of 
commutative algebras and let $M$ be an $A$-module. 
Then the associated chain complex of the simplicial $\Ga$-module
$C_*({\cal L}(X_*,M))\to {\cal L}(A,M)$ is a $\Y-$resolution.
\end{Le}

\S
\pn{\it Proof}. Since $\ee$ is weak equivalence of simplicial
algebras it is  a homotopy equivalence in the category of
 simplicial vector spaces. 
Thus $M\t D^{\ll}(X_*)\to M\t D^{\ll}(A_*)$ is also a homotopy
equivalence, for any particion $\ll$. It follows that
$${\cal L}(X_*,M)([n])^{\Sigma (\ll)}\to {\cal L}(A,M)([n])^{\Sigma (\ll)}$$
is also a homotopy equivalence of simplicial vector spaces. \qed

\B
The following is our main result.

\Thm{ For any commutative ring $A$ and any $A$-module $M$, there is 
a canonical isomorphism
$${\sf D}_i(A,M)\cong \pi_i^{\Y}({\cal L}(A,M)),\ i\geq 0$$
between the Andr\'e-Quillen homology and relative derived functors of
$\pi_0$ applied on the functor ${\cal L}(A,M)$.}

\S
\pn {\it Proof}. Thanks to Lemma \ref{nuli} 
the result is true for $i=0$. First consider the case,
when $M=A=K[t]$.
In this case Andr\'e-Quillen homology vanishes in positive dimensions
 by definition. On the other hand 
${\cal L}(K[t],K[t])$ is $\Y$-projective thanks to
Example \ref{patara} and 
therefore  $\pi_i^{\Y}({\cal L}(A,M))$ vanishes for all $i>0$.
One can use Lemma \ref{namravli} to 
conclude that $\pi_i^{\Y}({\cal L}(A,A))$ vanishes 
for all $i>0$ provided $A$ is a polynomial algebra. For the 
next step, we proof 
that the result is true if $A$ is a polynomial algebra and $M$ is
any $A$-module. 
We have to prove that  
$\pi_i^{\Y}({\cal L}(A,M))$ also vanishes
for $i>0$. We already proved this fact if $M=A$.
By additivity the functor $\pi_i^{\Y}({\cal L}(A,-))$
vanishes on free $A$-modules. By Lemma \ref{sizuste} the functor 
$\pi_*^{\Y}({\cal L}(A,-))$ assignes the long exact sequence to a
short exact sequence of $A$-modules. Therefore 
we can consider such an exact sequence
associated to a short exact sequence of $A$-modules
$$0\to N\to F\to M\to 0$$
with free $F$. Since the result is true if $i=0$ one obtains by
induction on $i$, that $\pi_i^{\Y}({\cal L}(A,M))=0$ provided $i>0$.
Now consider the  general case. Let $P_*\to A$ be a free simplicial 
resolution in the category of commutative algebras. Then we have
$$\Omega ^1_{P_*}\t _{P_*}M\cong \pi_0^{\Y}({\cal L}(P_*,M))$$
Thanks to Lemma \ref{rezolventa} 
$C_*({\cal L}(P_*,M))\to {\cal L}(A,M)$ is a $\Y$-resolution
consisting with $\pi_*^{\Y}$-acyclic objects and the result follows.
 \qed

The main theorem together with the main result of 
\cite{piri} yields:
\Cor{If $Char(K)=0$, then for any commutative  algebra $A$ and any 
$A$-module $M$ one has a natural isomorphism
$${\sf D}_*(A,M)\cong {\sf H}^{\G}_{*}(A,M).$$
}This  
fact was also proved  in  \cite{Sarah} 
based on the combinatorical and homotopical  
analysis of the space of fully grown trees.

\B
\pn {\bf Remarks}. i) We let $t:\Ga^{op}\to {\sf Vect}$ be the
functor which assignes the vector space of all maps $f:[n]\to K$, $f(0)=0$
to $[n]$. Then $t\t_{\G}F\cong \pi_0(F)$ (see Proposition 2.2
 \cite{hodge}). Hence $\pi_*^{\Y}$ van be also defined as the relative derived
functors of the functor $t\t_{\G}(-):\Ga$-$nod\to \sf Vect$. More generally
one can take any functor $T:\Ga^{op}\to {\sf Vect}$ and define
${\sf Tor}^{\Y}_*(T,F)$ as the value of the relative derived functors
(with respect of $\Y$-exact sequences) of the functor 
$T\t_{\G}(-):\Ga$-$nod\to \sf Vect$. Then our result claims
that
$${\sf D}_*(A,M)\cong {\sf Tor}^{\Y}_*(t, {\cal L}(A,M)).$$
Based on the Proposition 1.15 of \cite{hodge} the argument given above
shows that
$${\sf D}_*^{\{n\}}(A,M)\cong {\sf Tor}^{\Y}_*(\Lambda ^n \circ t, {\cal L}(A,M)),$$
where ${\sf D}_*^{\{n\}}(A,M)$ are defined using K\"ahler $n$-differentials:
$${\sf D}_*^{\{n\}}(A,M):= H_*(C_*(\Omega ^n_{P_*}\t _{P_*}M))$$
and for $n=1$ one recovers the main theorem.

\S
ii) All results remains true if $K$ is any commutative ring
and $A$ and $M$ are projective as $K$-mudule.


\bigskip

\centerline{\bf Acknowledgments}

\bigskip
\pn This work was written during my visit at MPI at Bonn.

\pn Razmadze Mathematical Institute \pn
Rukhadze str. 1 Tbilisi 380093 \pn 
Republic of Georgia \pn
{\tt email: pira@rmi.acnet.ge}


\begin{thebibliography}{999999999}
\bigskip

\bibitem{BF} {\sc A. K. Bousfield} and {\sc E.M. Friedlander}. {\it Homotopy  
theory of $\Gamma $-spaces, spectra, and bisimplicial sets}. Geometric
applications of homotopy theory (Proc. Conf., Evanston, Ill., 1977),  
II, pp. 80--130, Lecture Notes in Math., 658, Springer, 1978.

\M
\bibitem {DP} {\sc A. Dold} and {\sc D. Puppe}. {\it Homologie nicht-additiver 
Funktoren, Anwendungen.}  Ann. Inst. Fourier {\bf 11} (1961), 201--312.




\M
\bibitem{JLL} {\sc J. - L. Loday}. {\it Op\'erations sur l'homologie  
cyclique des $alg\grave ebres$ commutatives. } Invent. math. {\bf  
96} (1989), 205--230.


\M
\bibitem{homology} {\sc S. Mac Lane}. Homology. Classics in Mathematics. 
Springer-Verlag, Berlin, 1995. x+422 pp

\M
\bibitem{hodge} {\sc T. Pirashvili}. {\it Hodge decomposition of higher  
order Hochschild homology}.  Ann. Sci. \'Ecole Norm. Sup. (4) 33 (2000), no. 2, 151--179. 

\M
\bibitem{dold} {\sc T. Pirashvili}. {\it Dold-Kan type 
theorem for $\Gamma$-groups}. Math. Ann. 318 (2000), no. 2, 277--298. 


\M
\bibitem {piri} {\sc T. Pirashvili} and {\sc 
B. Richter.} {\it Robinson-Whitehouse complex and stable
homotopy}. Topology 39 (2000), no. 3, 525--530. 




\M
\bibitem{Q} {\sc D.G. Quillen.} {\it On the (co)homology of commutative rings.} 
 AMS Proc. Sym. Pure Math. {\bf XVII} (1970), 65--87.


\M


\M
\bibitem{segal} {\sc G. Segal}. {\it  Categories and cohomology theories.}  
Topology {\bf 13} (1974), 293--312.


\M
\bibitem{Sarah} {\sc  S. A. Whitehouse}. {\it Gamma (co)homology of commutative  
algebras and some related representations of the symmetric group}.  
Thesis. University of Warwick. 1994.




\end{thebibliography}
\end{document}